\documentclass[12pt]{article}
\usepackage{graphicx}
\usepackage{amsmath,amsthm,amssymb,enumerate}
\usepackage{euscript,mathrsfs}
\usepackage[left=3cm,right=3cm,top=3.5cm,bottom=3.5cm]{geometry}
\usepackage{color}
\catcode`\@=11 \@addtoreset{equation}{section}

\catcode`\@=12

\allowdisplaybreaks

\newtheorem{Theorem}{Theorem}[section]
\newtheorem{Proposition}{Theorem}[section]
\newtheorem{Lemma}{Theorem}[section]
\newtheorem{Corollary}{Theorem}[section]
\newtheorem{Definition}[Theorem]{Definition}

\newtheorem{Remark}[Theorem]{Remark}

\newcommand{\bTheorem}[1]{
\begin{Theorem} \label{T#1} }
\newcommand{\eT}{\end{Theorem}}

\newcommand{\bProposition}[1]{
\begin{Proposition} \label{P#1}}
\newcommand{\eP}{\end{Proposition}}

\newcommand{\bLemma}[1]{
\begin{Lemma} \label{L#1} }
\newcommand{\eL}{\end{Lemma}}

\newcommand{\bCorollary}[1]{
\begin{Corollary} \label{C#1} }
\newcommand{\eC}{\end{Corollary}}

\newcommand{\bRemark}[1]{
\begin{Remark} \label{R#1} }
\newcommand{\eR}{\end{Remark}}

\newcommand{\bDefinition}[1]{
\begin{Definition} \label{D#1} }
\newcommand{\eD}{\end{Definition}}

\newcommand{\cD}{c({\rm data})}

\newcommand{\vTe}{\Theta_{\ep}}

\newcommand{\bFormula}[1]{
\begin{equation} \label{#1}}
\newcommand{\eF}{\end{equation}}

\newcommand{\Ov}[1]{\overline{#1}}

\newcommand{\DC}{C^\infty_c}

\newcommand{\vr}{\varrho}
\newcommand{\vre}{\vr_\ep}

\newcommand{\vue}{\vu_\ep}

\newcommand{\vu}{\vc{u}}
\newcommand{\vc}[1]{{\bf #1}}

\newcommand{\Div}{{\rm div}_x}
\newcommand{\Grad}{\nabla_x}

\newcommand{\tn}[1]{\mathbb{#1}}
\newcommand{\dx}{\,{\rm d} {x}}

\newcommand{\dt}{\,{\rm d} t }

\newcommand{\intO}[1]{\int_{\Omega} #1 \ \dx}

\newcommand{\intR}[1]{\int_{R^3} #1 \ \dx}

\newcommand{\intTO}[1]{\int_0^T \int_{R^3} #1 \ \dx \dt}

\newcommand{\ep}{\varepsilon}

\definecolor{Cgrey}{rgb}{0.85,0.85,0.85}
\definecolor{Cblue}{rgb}{0.50,0.85,0.85}
\definecolor{Cred}{rgb}{1,0,0}
\definecolor{fancy}{rgb}{0.10,0.85,0.10}

\newcommand\Cbox[2]{%
    \newbox\contentbox%
    \newbox\bkgdbox%
    \setbox\contentbox\hbox to \hsize{%
        \vtop{
            \kern\columnsep
            \hbox to \hsize{%
                \kern\columnsep%
                \advance\hsize by -2\columnsep%
                \setlength{\textwidth}{\hsize}%
                \vbox{
                    \parskip=\baselineskip
                    \parindent=0bp
                    #2
                }%
                \kern\columnsep%
            }%
            \kern\columnsep%
        }%
    }%
    \setbox\bkgdbox\vbox{
        \color{#1}
        \hrule width  \wd\contentbox %
               height \ht\contentbox %
               depth  \dp\contentbox
        \color{black}
    }%
    \wd\bkgdbox=0bp%
    \vbox{\hbox to \hsize{\box\bkgdbox\box\contentbox}}%
    \vskip\baselineskip%
}

\date{}

\begin{document}


\title{An anelastic approximation arising in astrophysics}

\author{Donatella Donatelli
\thanks{The research of D.D. was partially supported by European Union MSCA-ITN-2014-ETN - Marie Sklodowska-Curie Innovative
Training Networks (ITN-ETN) Project: grant agreement 642768 ModCompShock}  \and Eduard Feireisl
\thanks{The research of E.F. leading to these results has received
funding from the European Research Council under the European Union's Seventh Framework
Programme (FP7/2007-2013)/ ERC Grant Agreement 320078. The Institute of Mathematics of the Academy of Sciences of
the Czech Republic is supported by RVO:67985840. E.F. thanks Gran Sasso Science Institute for hospitality and support during
his stay in L'Aquila.}}


\maketitle

\bigskip

\centerline{Department of Information Engineering, Computer Science and Mathematics}

\centerline{University of L'Aquila, 67100 L'Aquila, Italy}

\bigskip

\centerline{Institute of Mathematics of the Academy of Sciences of the Czech Republic}

\centerline{\v Zitn\' a 25, CZ-115 67 Praha 1, Czech Republic}

\bigskip


\begin{abstract}

We identify the asymptotic limit of the compressible non-isentropic Navier-Stokes system in the regime of low Mach, low Froude and high
Reynolds number. The system is driven by a long range gravitational potential. We show convergence to an anelastic system for ill-prepared initial data. The proof is based on frequency localized Strichartz estimates for the acoustic equation based on the recent work of Metcalfe and Tataru.

\end{abstract}



\section{Introduction}
\label{i}

A large class of phenomena arising in science occurs under various scaling regimes, among them  it is particularly relevant  the Low Mach number regime where the characteristic velocity of the system  is much smaller than the sound speed. It becomes, then,   very important to understand the separation of scales that occurs under this specific state and to identify the asymptotic regime (``low Mach Number model''). In particular we can interpret  the final system as a set of equations supporting the equilibria  infinite fast oscillating waves.

In this paper we consider the low Mach number regime for supernovas, for the  model see Almgren et al. \cite{ABMZ}, \cite{ABRZ1}, \cite{ABRZ2}, where it is described the evolution of the so called Type Ia supernovae, but the same set of equation can also be applied to any other models in astrophysics. With very simple words we can say that a  supernova is the explosion of a star. These type of events are characterized  by a range in timescales, in fact  it takes 100 years for the convection process  that precedes ignition while the duration of the explosion is of one second and, in the last minutes of the convective phase, velocity reachs almost the  $1\%$ of the sound speed. For this reasons it is appropriate to talk about low Mach number regime when dealing with supernovae.

In this paper we consider a very simplified hydrodinamical model describing the above mentioned phenomena, we will not take into account reaction phenomena and we will consider very easy state equation for the pressure and the internal energy. Our \emph{primitive system} looks as follows 

\begin{eqnarray}
\label{i5} \partial_t \vr + \Div (\vr \vu) &=& 0 \\
\label{i6} \partial_t (\vr \vu) + \Div (\vr \vu \otimes \vu ) + \frac{1}{\ep^2} \Grad (\vr \Theta)^\gamma &=& \ep^\alpha \Div \mathbb{S} (\Grad \vu) + \frac{1}{\ep^2}
\vr \Grad F\\
\label{i7} \partial_t (\vr \Theta) + \Div (\vr \Theta \vu) &=& 0 \\
\label{i8} \mathbb{S} (\Grad \vu) = \left( \Grad \vu + \Grad^t \vu - \frac{2}{3} \Div \vu \mathbb{I} \right) &+& \lambda \Div \vu \mathbb{I},\
\lambda \geq 0,
\end{eqnarray}
supplemented with the far-field conditions
\begin{equation} \label{i9}
\vr \to \Ov{\vr}, \ \vu \to 0, \ \Theta \to 1 \ \mbox{as}\ |x| \to \infty,
\end{equation}
where the parameter $\ep \to 0$, and $\alpha > 0$.

Our goal is to identify the asymptotic limit of the system (\ref{i5}-\ref{i8}) as

\begin{eqnarray}
\Div (\vr_0 \vc{V}) &=& 0 \label{i1} \\
\partial_t \vc{V} + \vc{V} \cdot \Grad \vc{V} + \Grad \Pi &=& - \frac{\vr_0}{\mathcal{R}} \Grad F \label{i2} \\
\partial_t (\mathcal{R}) + \Div (\mathcal{R} \vc{V}) &=& 0 \label{i3}
\end{eqnarray}
with the unknown density $\mathcal{R}$, velocity $\vc{V}$, and the pressure $\Pi$, see Almgren et al. \cite{ABMZ}, \cite{ABRZ1}, \cite{ABRZ2}. The background density profile
$\vr_0$ is determined as the unique solution of the static equation
\begin{equation} \label{i4}
\Grad \vr_0^\gamma = \vr_0 \Grad F \ \mbox{in}\ R^3, \ \vr_0 \to \Ov{\vr} > 0 \ \mbox{as}\ |x| \to \infty.
\end{equation}

A similar system, written in terms of the potential temperature $\mathcal{T} = \vr_0 / \mathcal{R}$, has been identified
as a low Mach number approximation of strongly stratified fluid flows in meteorology, see \cite{FeKlNoZa}.
The principal difference between \cite{FeKlNoZa} and the present problem is the geometry of the physical space; an infinite strip is relevant in meteorology, while
the whole Euclidean space $R^3$ is considered in astrophysics mimicking a large neighborhood of a gaseous star.
Although this might seem as a minor modification, the related difficulties stemming from the action of acoustic waves in the 
incompressible limit must be handled in a rather different manner. Here, we use a seminal paper of Metcalfe and Tataru \cite{MetTat1} on wave propagation 
with variable wave speed. Combining their method with the local frequency estimates to control the acoustic waves propagation 
may be seen as the main contribution and the main added value of the present paper in comparison with
\cite{FeKlNoZa}.

In accordance with the previous discussion, the potential $F$ shares
the asymptotic properties with the standard gravitational potential:
\begin{eqnarray}
\label{i4a} F \in C^\infty(R^3),\ F(x) &>& 0 \ \mbox{for all}\ x \in R^3, \\
\label{i4b} \underline{F}\frac{1}{|x|} &\leq& F(x) \leq \Ov{F} \frac{1}{|x|} \ \mbox{for all}\ |x| > R,\\
\label{i4c} |x|^2 |\Grad F(x) | + |x|^3 |\Grad^2 F(x)| &\leq& c \ \mbox{for all}\ x \in R^3.
\end{eqnarray}

Equations (\ref{i5}--\ref{i7}) represent a scaled isentropic Navier-Stokes system written in terms of the density $\vr$, the velocity $\vu$, and the potential
temperature $\Theta$. The relevant value of the adiabatic exponent is $\gamma = \frac{5}{3}$ for the monoatomic gas state equation. The singular scaling
factors in (\ref{i6}) are the Mach number and the Froude number proportional to $\ep^2$ and the Reynolds number proportional to $\ep^{-\alpha}$. The relevance of the regime $\ep \to 0$ in astrophysical models is discussed in \cite{ABRZ1}.

We consider the singular limit $\ep \to 0$ for the \emph{ill-prepared} initial data:
\begin{eqnarray}
\label{i10} \vr(0, \cdot) &=& \vr_0 + \ep \vr^{(1)}_{0,\ep}, \ \| \vr^{(1)}_{0,\ep} \|_{L^1 \cap L^\infty (R^3)} \leq c,
\vr^{(1)}_{0, \ep} \to \vr^{(1)}_0 \ \mbox{in}\ L^1(R^3), \\
\label{i11} \vu(0, \cdot) &=& \vu_{0,\ep}, \ \| \vu_{0,\ep} \|_{L^2 \cap L^\infty (R^3)} \leq c, \ \vu_{0, \ep} \to \vc{v}_0 \
\mbox{in}\ L^2(R^3), \\
\label{i12} \Theta(0, \ep) &=& 1 + \ep^2 \Theta_{0, \ep}^{(2)}, \ \| \Theta^{(2)}_{0,\ep} \|_{L^1 \cap L^\infty (R^3)} \leq c, \
\Theta^{(2)}_{0, \ep} \to \Theta^{(2)}_0 \ \mbox{in}\ L^1(R^3).
\end{eqnarray}
In particular, we do not require the limit velocity field $\vc{v}_0$ to satisfy the anelastic constraint (\ref{i1}).

In contrast with \cite{FeKlNoZa}, where the problem for ill-prepared data was considered only in the viscous regime (constant positive Reynolds number),
we are able to handle the vanishing viscosity limit in the presence of the disturbing effect of acoustic waves generated by the ill-prepared data.
For this kind of problem, the lack
of compactness of the velocity must be compensated by the structural stability properties of the system encoded in the relative energy inequality,
identified for the present system in \cite{FeKlNoZa}.

As the underlying physical space is unbounded, the effect of acoustic waves is likely to be annihilated by dispersion. Unfortunately,
to our best knowledge, no direct dispersive estimates are available for the acoustic system in question.
The problem is that we need rather strong \emph{global in time and space} estimates of Strichartz's type for a non-constant coefficient wave operator
\[
\mathcal{A}_{\vr_0}: v \mapsto - \frac{p'(\vr_0)}{\vr_0} \Div \left( \vr_0 \Grad v \right),
\]
with $\vr_0$ generated by the long-range potential $F$. Luckily, however, it is enough to have these estimates for a ``dense'' set of initial data,
in particular, we may work on a frequency localized space, where the low and high frequencies are cut-off. Our approach,
reminiscent of the decomposition technique used by Smith and Sogge \cite{SmSo}, can be summarized as follows.
\begin{itemize}
\item Use the result of De Bi{\`e}vre and Pravica \cite{DeBpr1}, \cite{DeBpr2} to show that the point spectrum of $\mathcal{A}_{\vr_0}$ is empty and that
$\mathcal{A}_{\vr_0}$ satisfies the Limiting absorption principle on a weighted Hilbert space.
\item Show frequency and space localized decay estimates by means of the method proposed in \cite{EF80}. These are estimates on solutions of the
wave equation generated by $\mathcal{A}_{\vr_0}$ localized in both the physical and frequency space.
\item Combine the local estimates with the global result of Metcalfe and Tataru \cite{MetTat1}. As a result we obtain frequency localized Strichartz estimates
    for $\mathcal{A}_{\vr_0}$ that may be of independent interest, see Section \ref{A}.
\end{itemize}

The paper is organized as follows. We start with some preliminary material in Section \ref{m}, introducing the concept of
dissipative solution to the primitive system, discussing solvability of the target system, and, finally, stating the main result.
In Section \ref{r}, we introduce the relative energy inequality and derive the necessary uniform bounds
for solutions of the primitive system independent of the scaling parameters.
Section \ref{A} is the heart of the paper. We study the underlying acoustic equation and establish the dispersive estimates
of Strichartz's type for frequency localized data. The proof of the main result is
completed in Section \ref{C}.

\section{Preliminaries, main result}
\label{m}

We introduce the concepts of \emph{dissipative weak solution} for the primitive system (\ref{i5}--\ref{i9}) and
\emph{strong solution} for the target system (\ref{i1}--\ref{i4}).
Solutions of both systems will be considered on the whole physical space $R^3$, with relevant far field conditions for $|x| \to \infty$.

\subsection{Dissipative solutions of the primitive system}

We consider the weak solutions of system (\ref{i5}--\ref{i9}) in the class
\begin{eqnarray}
\label{m1}
( \vr - \Ov{\vr} ) &\in&  C_{\rm weak} ([0,T]; (L^2 + L^\gamma)(R^3)) \\
\label{m2}
\vr \vu &\in& C_{\rm weak} ([0,T]; (L^2 + L^{\frac{2\gamma}{\gamma + 1}})(R^3;R^3)) \\
\label{m3}
\vr( \Theta - 1) &\in& C_{\rm weak} ([0,T]; (L^2 + L^\gamma)(R^3)), \ \Theta \in L^\infty((0,T) \times R^3),
\end{eqnarray}
satisfying
\begin{equation} \label{m4}
\begin{split}
\left[ \intR{ \vr \varphi  } \right]_{t = \tau_1}^{t = \tau_2} &= \int_{\tau_1}^{\tau_2}
\intR{ \left[ \vr \partial_t \varphi + \vr \vu \cdot \Grad \varphi \right] }  \dt\\
\mbox{for any}\ 0 \leq \tau_1 \leq \tau_2 \leq T, \ \mbox{and for any} \ &\varphi \in C^1_c([0,T] \times R^3);
\end{split}
\end{equation}
\begin{equation} \label{m5}
\begin{split}
\left[ \intR{ \vr \vu \cdot \varphi  } \right]_{t = \tau_1}^{t = \tau_2} &= \int_{\tau_1}^{\tau_2}
\intR{ \left[ \vr \vu \cdot \partial_t \varphi + \vr \vu \otimes \vu : \Grad \varphi + \frac{1}{\ep^2} (\vr \Theta)^\gamma \Div \varphi \right] }  \dt\\
&- \int_{\tau_1}^{\tau_2} \intR{ \left[ \ep^\alpha \tn{S} (\Grad \vu) : \Grad \varphi - \frac{1}{\ep^2} \vr \Grad F \cdot \varphi \right] } \ \dt \\
\mbox{for any}\ 0 \leq \tau_1 &\leq \tau_2 \leq T,\ \mbox{and for any} \ \varphi \in C^1_c([0,T] \times R^3; R^3);
\end{split}
\end{equation}
\begin{equation} \label{m6}
\begin{split}
\left[ \intR{ \vr G(\Theta) \varphi  } \right]_{t = \tau_1}^{t = \tau_2} &= \int_{\tau_1}^{\tau_2}
\intR{ \left[ \vr G(\Theta) \partial_t \varphi + \vr G(\Theta) \vu \cdot \Grad \varphi \right] }  \dt\\
\mbox{for any}\ 0 \leq \tau_1 &\leq \tau_2 \leq T,\ \mbox{and for any} \ \varphi \in C^1_c([0,T] \times R^3) \ \mbox{and any}\ G \in C(R).
\end{split}
\end{equation}

\begin{Remark} \label{Rm1}
Equation (\ref{m6}), introduced in \cite{FeKlNoZa}, can be viewed as a renormalized formulation of (\ref{i7}).
\end{Remark}

In this paper, a slightly different form of renormalization of (\ref{i7}) is needed, namely
\begin{equation} \label{m6a}
\begin{split}
\left[ \intR{ b(\vr \Theta) \varphi  } \right]_{t = \tau_1}^{t = \tau_2} &= \int_{\tau_1}^{\tau_2}
\intR{ \left[ b(\vr \Theta) \partial_t \varphi + b(\vr  \Theta) \vu \cdot \Grad \varphi \right] }  \dt\\
&+ \int_{\tau_1}^{\tau_2}
\intR{ \left[ b(\vr \Theta) - b'(\vr \Theta) (\vr \Theta) \right] \Div \vu } \ \dt
\\
\mbox{for any} \ &\varphi \in \DC([0,T] \times R^3) \ \mbox{and any}\ b \in C^1(R), \ b' \in C_c(R).
\end{split}
\end{equation}

In addition, the
\emph{dissipative solutions} are characterized by the energy inequality
\begin{equation} \label{m7}
\begin{split}
&\left[ \intR{ \left[ \frac{1}{2} \vr |\vu|^2 + \frac{H (\vr \Theta) - H' (\vr_0)(\vr - \vr_0) - H(\vr_0) }{\ep^2} \right]  } \right]_{t = 0}^{t = \tau} \\
& + \ep^\alpha \int_0^\tau \intR{ \mathbb{S}(\Grad \vu) : \Grad \vu } \ \dt \leq 0  \\
\mbox{for a.a.} \ \tau &\in (0,T), \ \mbox{where}\ H (Z) = \frac{1}{\gamma - 1}Z^\gamma.
\end{split}
\end{equation}

The \emph{existence} of dissipative weak solution to problem (\ref{i5}--\ref{i9}) on a bounded domain with the non-slip
boundary conditions was shown by Maltese et al. \cite{MMMNPZ}.

\subsection{Smooth solutions of the target problem}

Our technique requires the target system to admit a smooth solution. Although
there is probably no precise statement concerning solvability of the target system (\ref{i1}--\ref{i3}) available in the literature,
we may anticipate, also in view of
the work by Oliver \cite{Oli} on the anelastic Euler system, the existence of local-in-time smooth solutions in the class:
\begin{eqnarray}
\label{m8}
\vc{V} &\in& C([0,T_{\rm max}); W^{m,2}(R^3;R^3)) \\
\label{m9}
\Pi &\in& C([0,T_{\rm max}); W^{m,2}(R^3;R^3)) \\
\label{m10}
\mathcal{R} &\in& C([0,T_{\rm max}); W^{m,2}(R^3))
\end{eqnarray}
for $m > 3$, defined on a maximal time interval $[0, T_{\rm max})$.

\subsection{Main result}

In order to state our main result we need a weighted variant of the Helmholtz decomposition,
\[
\vc{v} = \vc{H}_{\vr_0}[\vc{v}] + \Grad \Phi,
\]
where the potential $\Phi \in D^{1,2}(R^3)$ is uniquely determined by
\[
\Div \left( \vr_0 \Grad \Phi \right) = \Div ( \vr_0 \vc{v}), \
\Phi \in D^{1,2}(R^3) \ \mbox{for}\ \vc{v} \in L^2(R^3,R^3).
\]

Our goal is to show the following result.

\begin{Theorem} \label{Tm1}
Let the potential $F$ satisfy (\ref{i4a}--\ref{i4c}). Let $\Ov{\vr} > 0$ and let $\vr_0$ be the unique solution of (\ref{i4}). In addition, suppose
\begin{equation} \label{HYP}
\gamma > \frac{3}{2}, \  0 < \alpha < \frac{4}{3}.
\end{equation}
Let $\{ \vre, \vue, \vTe \}_{\ep > 0}$ be a family of dissipative solutions to the primitive system (\ref{i5}--\ref{i9}) in $(0,T) \times R^3$, with the
initial values satisfying (\ref{i10}--\ref{i12}). Let the target system (\ref{i1}--\ref{i3}) admit a strong solution $[\vc{V}, \Pi,
\mathcal{R}]$ in $(0,T) \times R^3$ in the class
(\ref{m8}--\ref{m10}), with the initial values
\[
\vc{V}(0, \cdot) = \vc{H}_{\vr_0} [\vc{v}_0],\ \mathcal{R}(0, \cdot) = \frac{\vr_0}{\Theta^{(2)}_0 } > 0.
\]

Then
\begin{equation} \label{m11}
\begin{split}
\sup_{t \in [0,T]} \| \vr_\ep (t, \cdot) - \vr_0 \|_{(L^2 + L^\gamma)(R^3)} &\to 0 \\
\sup_{t \in [0,T]} \left\| \vTe (t, \cdot) - \frac{\vr_0}{ \mathcal{R} }(t, \cdot) \right\|_{L^2(R^3)} &\to 0, \\
\int_0^T \left\| \sqrt{\frac{\vre}{\vr_0}} \vue -   \vc{V} \right\|^2_{L^2(K)} \dt &\to 0 \ \mbox{for any compact}\ K \in R^3
\end{split}
\end{equation}
as $\ep \to 0$.
\end{Theorem}

\begin{Remark} \label{RM2}

Note that the result is \emph{path dependent}, meaning the values of the Reynolds and Mach/Froude numbers are interrelated. In particular, the Reynolds number
cannot be too large with respect to the Mach number.

\end{Remark}

The rest of the paper is devoted to the proof of Theorem \ref{Tm1}.

\section{Relative energy, uniform bounds}
\label{r}

Similarly to \cite{FeKlNoZa}, we introduce the \emph{relative energy} functional
\[
\begin{split}
\mathcal{E} \left( \vr, \Theta, \vu \ \Big| \ r , \vc{U} \right) =
\intR{  \left[ \frac{1}{2} \vr |\vu - \vc{U}|^2 + \frac{ H(\vr \Theta) - H'(r)(\vr \Theta - r) - H(r)}{\ep^2} \right]}  \\
\end{split}
\]
for any pair of ``test functions''
\begin{equation} \label{class}
(r - \Ov{\vr}) \in \DC([0,T] \times R^3), \ r > 0,\ \vc{U} \in \DC([0,T] \times R^3; R^3).
\end{equation}

\subsection{Relative energy inequality and energy estimates}

As shown in \cite{FeKlNoZa} (cf. also \cite{FeJiNo1}) any dissipative solution to the primitive system (\ref{m4}--\ref{m7}) satisfies the
\emph{relative energy inequality}
\begin{equation} \label{r1}
\begin{split}
\left[ \mathcal{E} \left( \vr, \Theta, \vu \Big| r, \vc{U} \right) \right]_{t = 0}^{t = \tau}
&+ \ep^\alpha \int_{0}^{\tau}
\intR{ \mathbb{S} (\Grad (\vu-\vc{U})) : \Grad (\vu - \vc{U})} \dt
\\
\leq
\int_{0}^{\tau} &
\intR{  \vr \left( \partial_t \vc{U} + \vu \cdot \Grad \vc{U} \right) \cdot \left( \vc{U} - \vu \right)
+ \ep^\alpha \mathbb {S} (\Grad \vc{U}) : \Grad (\vc{U} - \vu ) } \dt
\\
+ \frac{1}{\ep^2} \int_0^\tau &  \intR{
\left[ (r - \vr \Theta) \partial_t H'(r) + \Grad H'(r) \cdot (r \vc{U} - \vr \Theta \vu ) \right] } \dt
\\
-  \frac{1}{\ep^2} \int_0^\tau & \intR{ \left[ \Div \vc{U} \Big( (\vr \Theta)^\gamma - r^\gamma \Big)
     + \vr \Grad F \cdot (\vc{U} - \vu) \right] } \dt
\end{split}
\end{equation}
for all sufficiently smooth functions belonging to the class (\ref{class}).

\begin{Remark} \label{Rclass}

Here, similarly to (\cite{FeKlNoZa}), the class of admissible test functions (\ref{class}) can be considerably extended by means of a density argument.

\end{Remark}

The relative energy inequality (\ref{r1}), together with the renormalized equation (\ref{m6}), give rise to a number of uniform bounds for the
family $\{ \vre, \vTe, \vue \}_{\ep > 0}$ independent of the scaling parameter. First, let us introduce some notation. Each
measurable function $h$ can be decomposed as $h=[h]_{{\rm ess}}+[h]_{{\rm res}}$, where
\[
\begin{split}
[h]_{{\rm ess}}&=\chi(\vr_\ep\Theta_\ep)h,\ [h]_{{\rm res}}=(1-\chi(\vre\Theta_\ep))h\\
\mbox{where}\ &\chi \in \DC (0, \infty),\ \chi \geq 0,\
\ \chi (Y) = 1\quad \mbox{whenever}\  \frac{1}{2} (\min_{x \in R^3} \varrho_0(x)) \leq Y \leq 2 (\max_{x \in R^3} \vr_0(x)).
\end{split}
\]
Under the hypotheses (\ref{i10}--\ref{i12}) imposed on the initial data,
the following estimates were proved in \cite[Section 3.2]{FeKlNoZa}:
\begin{eqnarray}
\label{r2}
{\rm ess} \sup_{t \in [0,T]} \| \sqrt{\vre} \vue \|_{L^2(R^3; R^3)} &\leq& \cD,
\\
\label{r3}
\ep^{\alpha/2} \left\| \Grad \vue + \Grad^t \vue - \frac{2}{3} \Div \vue \tn{I} \right\|_{L^2((0,T) \times R^3; R^{3 \times 3})}
&\leq& \cD,\\
\label{r4}
\ep^{\alpha/2}  \left\| \Div \vue \right\|_{L^2((0,T) \times R^3)}
&\leq& \cD,\\
\label{r5}
{\rm ess} \sup_{t \in [0,T]} \left( \left\| \frac{ \vTe - 1 }{\ep^2} \right\|_{L^1(R^3)} + \left\| \frac{ \vTe - 1 }{\ep^2} \right\|_{L^\infty(R^3)} \right) &\leq& \cD \\
\label{r6}
{\rm ess} \sup_{t \in [0,T]} \left( \left\| \left[ \frac{ \vre - \vr_0 }{\ep} \right]_{\rm ess} \right\|_{L^2(R^3)} + \left\|
\left[ \frac{ \vre \vTe - \vr_0 }{\ep} \right]_{\rm ess} \right\|_{L^2(R^3)} \right) &\leq& \cD\\
\label{r7}
{\rm ess} \sup_{t \in [0,T]} \intR{ \left( [1]_{\rm res} + \left| [\vre ]_{\rm res} \right|^\gamma + \left| [\vre \vTe ]_{\rm res} \right|^\gamma  \right) }
&\leq& \ep^2 \cD.
\end{eqnarray}

\begin{Remark} \label{rR1}
In (\ref{r5}) and whenever convenient, we set $\vTe = 1$ on the vacuum set $\{ \vre = 0 \}$.
\end{Remark}

Finally, we use (\ref{r2}), (\ref{r3}), (\ref{r6}) to conclude
\begin{equation}
\label{r8}
\left\| \vue \right\|_{L^2(0,T; W^{1,2} (R^3; R^3))} \leq \ep^{-\alpha/2} \cD.
\end{equation}

\subsection{Pressure estimates}

The pressure estimates are one of the crucial ingredients of the existence theory for the compressible Navier-Stokes system, see Lions \cite{LI4}.
Similarly to Masmoudi \cite{MAS2} (cf. also \cite{FeJiNo1}), they are also needed to control some terms in the anelastic limit.
We use the quantities
\[
\begin{split}
\varphi (t,x) &= \phi(x) \Grad \Delta^{-1}_x \left[ b \left( \vre \vTe \right) \right]_{\rm res} , \ \phi \in \DC(R^3), \ \phi \geq 0 \\
b \in C^\infty[0, \infty),\ b &\geq 0, \ b(Y) = a Y^\beta, \ a \geq 0,\ 0 \leq \beta < \gamma \ \mbox{for}\ Y >>1
\end{split}
\]
as test functions in the momentum equation (\ref{m5}). Here, the inverse $\Delta^{-1}_x$ is defined in the standard way
\[
\Delta^{-1}_x [v] (x) = \frac{1}{4 \pi} \int_{R^3} \frac{v(y) }{|x - y |} \ {\rm d}y.
\]
Note that (\ref{r7}) yields
\begin{equation} \label{r9}
{\rm ess} \sup_{t \in (0,T)}
\left\| \left[ b \left( \vre \vTe \right) \right]_{\rm res} \right\|_{L^{q}(R^3)} \leq c \ep^{\frac{2}{q}} \
\mbox{for any} \ 1\leq q \leq \frac{\gamma}{\beta}.
\end{equation}

Following \cite[Section 4.2]{FeJiNo1} and denoting
\[
p(Z) = Z^\gamma,
\]
we deduce the relation
\begin{equation} \label{r10}
\frac{1}{\ep^2} \intTO{ \phi \Big( p(\vre \vTe) - p(\vr_0) \Big) \left[ b \left( \vre \vTe \right) \right]_{\rm res} } = \sum_{j=1}^7 I_{j,\ep} ,
\end{equation}
with
\[
\begin{split}
I_{1,\ep} &= \frac{1}{\ep^2} \intTO{ \Big( p(\vr_0) - p(\vre \vTe) \Big) \Grad \phi \cdot \Grad \Delta^{-1}_x \left[ b \left( \vre \vTe \right) \right]_{\rm res} },
\\
I_{2,\ep} &= - \frac{1}{\ep^2} \intTO{ \phi (\vre - \vr_0) \Grad F \cdot \Grad  \Delta^{-1}_x \left[ b \left( \vre \vTe \right) \right]_{\rm res} }
\\
I_{3,\ep} &= \ep^\alpha \intTO{ \mathbb{S} (\Grad \vue) : \Grad \left( \phi \Grad \Delta^{-1}_x  \left[ b \left( \vre \vTe \right) \right]_{\rm res} \right) }
\\
I_{4,\ep} &= - \intTO{ \vre \vue \otimes \vue :
\Grad \left( \phi \Grad \Delta^{-1}_x \left[ b \left( \vre \vTe \right) \right]_{\rm res} \right) }
\\
I_{5,\ep} &= \intTO{ \phi \vre \vue \cdot  \Grad   \Delta^{-1}_x \Big[\Div( \left[ b \left( \vre \vTe \right) \right]_{\rm res} \vue) \Big]}
\\
I_{6,\ep} &= \intTO{ \phi \vre \vue \cdot
\Grad \Delta^{-1}_x  \Big[ \Big( \left[ b \left( \vre \vTe \right) \right]_{\rm res}  - \left[ b \left( \vre \vTe \right) \right]_{\rm res} ' \vre \vTe \Big) \Div \vue \Big]   }
\\
I_{7,\ep} &= \left[ \intO{ \phi \vre \vue \cdot \Grad \Delta^{-1}_x \left[ b \left( \vre \vTe \right) \right]_{\rm res} (\tau, \cdot) } \right]_{\tau = 0}^{\tau = T}
\end{split}
\]
where we have also used the renormalized equation (\ref{m6a}).

We consider two complementary choices of $b$, namely
\begin{equation} \label{r11}
b(Y) = 0 \ \mbox{if}\ 0 \leq Y \leq \min_{x \in R^3} \vr_0, \ \mbox{or}\  b(Y) = 0 \ \mbox{if}\ Y \geq \max_{x \in R^3} \vr_0.
\end{equation}
Note that in both cases the integrand on the left-hand side of (\ref{r10}) has a definite sign. Now we proceed in several steps.

\subsubsection{Integrals $I_1-I_3$, $I_7$}

We focus on the (more difficult) former case in (\ref{r11}) obtaining
\[
|I_{1,\ep}| \leq \frac{c}{\ep^2} \intTO{ \left|  \Big[ p(\vr_0) - p(\vre \vTe) \Big]_{\rm ess} \right| | \Grad \Delta^{-1}_x[b(\vre \vTe)]_{\rm res}
| |\Grad \phi | }
\]
\[
+ \frac{c}{\ep^2} \intTO{ \left|  \Big[ p(\vr_0) - p(\vre \vTe) \Big]_{\rm res} \right| | \Grad \Delta^{-1}_x[b(\vre \vTe)]_{\rm res}||\Grad \phi | }.
\]
By virtue of the Sobolev embedding relations,
\[
\begin{split}
\left\| \Grad \Delta^{-1}_x[b(\vre \vTe)]_{\rm res} \right\|_{L^2(R^3; R^3)} &\leq c \left\| [b(\vre \vTe)]_{\rm res} \right\|_{L^{6/5}(R^3)},\\
\left\| \Grad \Delta^{-1}_x[b(\vre \vTe)]_{\rm res} \right\|_{L^\infty(R^3; R^3)} &\leq \left\| [b(\vre \vTe)]_{\rm res} \right\|_{(L^{1} \cap L^q) (R^3)},
\ \mbox{for}\ q > 3.
\end{split}
\]
Thus for $0 < \beta < \frac{\gamma}{3}$, we may use (\ref{r9}), together with the uniform bounds (\ref{r6}), (\ref{r7}), to conclude that
\[
| I_{1,\ep} | \leq \ep^\omega \cD \ \mbox{for some}\ \omega > 0.
\]
Furthermore, in view of (\ref{r3}), (\ref{r6}), and (\ref{r7}), a similar argument can be applied to $I_{2,\ep}$, $I_{3, \ep}$, and $I_{7,\ep}$
to obtain
\[
| I_{j,\ep} | \leq \ep^\omega \cD \ \mbox{for some}\ \omega > 0, \ j=2,3,7.
\]

\subsubsection{Integrals $I_4$, $I_5$}

To handle $I_4$ we have to use the ``negative estimate'' (\ref{r8}) along with the embedding $W^{1,2}(R^3) \hookrightarrow L^6(R^3)$. Consequently,
we have
\[
\begin{split}
I_{4,\ep} &= - \intTO{ [ \vre \vue ]_{\rm ess} \otimes \vue :
\Grad \left( \phi \Grad \Delta^{-1}_x \left[ b \left( \vre \vTe \right) \right]_{\rm res} \right) }\\
& - \intTO{ [ \vre \vue ]_{\rm res} \otimes \vue :
\Grad \left( \phi \Grad \Delta^{-1}_x \left[ b \left( \vre \vTe \right) \right]_{\rm res} \right) },
\end{split}
\]
where, furthermore,
\[
\begin{split}
&\left| \intR{ [ \vre \vue ]_{\rm ess} \otimes \vue :
\Grad \left( \phi \Grad \Delta^{-1}_x \left[ b \left( \vre \vTe \right) \right]_{\rm res} \right)} \right|\\
& \leq c \| [ \vre \vue ]_{\rm ess} \|_{L^2(R^3;R^3)} \| \vue \|_{L^6(R^3; R^3)} \left\| \left[ b \left( \vre \vTe \right) \right]_{\rm res}
\right\|_{L^3(R^3)},
\end{split}
\]
and where, by virtue of (\ref{r2}), (\ref{r5}), (\ref{r8}), and (\ref{r9}), the right-hand side is of order $\ep^\omega$, $\omega > 0$ as $\ep^{2/3} \ep^{-\alpha /2} = \ep^\omega$, $\omega = \frac{2}{3} - \frac{\alpha}{2} > 0$ in accordance with hypothesis
(\ref{HYP}).

As for the residual part, we get
\[
\begin{split}
&\left| \intTO{ [ \vre \vue ]_{\rm res} \otimes \vue :
\Grad \left( \phi \Grad \Delta^{-1}_x \left[ b \left( \vre \vTe \right) \right]_{\rm res} \right) } \right|\\
&\leq \left\| [ \sqrt{\vre} ]_{\rm res} \right\|_{L^{2\gamma}(R^3)} \| \sqrt{\vre} \vue \|_{L^2(R^3; R^3)}
\| \vue \|_{L^6(R^3; R^3)} \left\| \left[ b \left( \vre \vTe \right) \right]_{\rm res} \right\|_{L^q(R^3)},\ q = \frac{6\gamma}{2\gamma - 3};
\end{split}
\]
Thus, in accordance with the bounds (\ref{r7}), (\ref{r8}), and (\ref{r9}), we need
\[
\frac{ 2 \gamma - 3}{3\gamma} + \frac{1}{\gamma} > \frac{\alpha}{2}, \ \mbox{meaning}, \ 0 < \alpha < \frac{4}{3},
\]
exactly as required by (\ref{HYP}).

The integral $I_{5,\ep}$ can be treated in the same manner.

\subsubsection{Integral $I_6$}

Finally,
\[
\begin{split}
|I_{6,\ep}| &\leq \left| \intTO{ \phi [\vre \vue]_{\rm ess} \cdot
\Grad \Delta^{-1}_x  \Big[ \Big( \left[ b \left( \vre \vTe \right) \right]_{\rm res}  - \left[ b \left( \vre \vTe \right) \right]_{\rm res} ' \vre \vTe \Big) \Div \vue \Big]   } \right|\\
&+ \left| \intTO{ \phi [\vre \vue]_{\rm res} \cdot
\Grad \Delta^{-1}_x  \Big[ \Big( \left[ b \left( \vre \vTe \right) \right]_{\rm res}  - \left[ b \left( \vre \vTe \right) \right]_{\rm res} ' \vre \vTe \Big) \Div \vue \Big]   } \right|,
\end{split}
\]
where, by means of the Sobolev embedding $W^{1,6/5} \hookrightarrow L^2$, we get
\[
\begin{split}
&\left| \intTO{ \phi [\vre \vue]_{\rm ess} \cdot
\Grad \Delta^{-1}_x  \Big[ \Big( \left[ b \left( \vre \vTe \right) \right]_{\rm res}  - \left[ b \left( \vre \vTe \right) \right]_{\rm res} ' \vre \vTe \Big) \Div \vue \Big]   } \right|\\
&\leq \left\| [\vre \vue]_{\rm ess} \right\|_{L^2(R^3;R^3)} \| \Grad \vu \|_{L^2(R^3; R^{3 \times 3})}
\left\|\left[ b \left( \vre \vTe \right) \right]_{\rm res} -  \left[ b \left( \vre \vTe \right) \right]_{\rm res} ' \vre \vTe \right\|_{L^3(R^3)}
\end{split}
\]
where the right-hand side can be estimated exactly is its counterpart in $I_{4, \ep}$.

To conclude
\[
\begin{split}
&\left| \intTO{ \phi [\vre \vue]_{\rm res} \cdot
\Grad \Delta^{-1}_x  \Big[ \Big( \left[ b \left( \vre \vTe \right) \right]_{\rm res}  - \left[ b \left( \vre \vTe \right) \right]_{\rm res} ' \vre \vTe \Big) \Div \vue \Big]   } \right|\\
&\leq \left\| [ \sqrt{\vre} ]_{\rm res} \right\|_{L^{2\gamma}(R^3)} \| \sqrt{\vre} \vue \|_{L^2(R^3; R^3)} \times \\
& \times \left\| \Grad \Delta^{-1}_x  \Big[ \Big( \left[ b \left( \vre \vTe \right) \right]_{\rm res}  - \left[ b \left( \vre \vTe \right) \right]_{\rm res} ' \vre \vTe \Big) \Div \vue \Big]   \right\|_{L^q(R^3; R^3)}, \ q = \frac{2 \gamma}{\gamma - 1},
\end{split}
\]
where, furthermore,
\[
\begin{split}
&\left\| \Grad \Delta^{-1}_x  \Big[ \Big( \left[ b \left( \vre \vTe \right) \right]_{\rm res}  - \left[ b \left( \vre \vTe \right) \right]_{\rm res} ' \vre \vTe \Big) \Div \vue \Big]   \right\|_{L^{2 \gamma/(\gamma - 1)} (R^3; R^3)}\\
& \leq \left\| \Big( \left[ b \left( \vre \vTe \right) \right]_{\rm res}  - \left[ b \left( \vre \vTe \right) \right]_{\rm res} ' \vre \vTe \Big) \Div \vue
\right\|_{L^{6\gamma/ (5\gamma - 3)}(R^3)},
\end{split}
\]
and, by H\" older's inequality,
\[
\begin{split}
&\left\| \Big( \left[ b \left( \vre \vTe \right) \right]_{\rm res}  - \left[ b \left( \vre \vTe \right) \right]_{\rm res} ' \vre \vTe \Big) \Div \vue
\right\|_{L^{6\gamma/ (5\gamma - 3)}(R^3)},
\\
&\leq \left\| \Div \vue \right\|_{L^2(R^3;R^3)} \left\| \left[ b \left( \vre \vTe \right) \right]_{\rm res}  - \left[ b \left( \vre \vTe \right) \right]_{\rm res} ' \vre \vTe  \right\|_{L^{6 \gamma/ (2\gamma - 3)}(R^3)};
\end{split}
\]
whence the required bound follows again from hypothesis (\ref{HYP}).

\subsubsection{Pressure estimates - conclusion}

We conclude that, under hypothesis (\ref{HYP}),
\begin{equation} \label{r12}
\int_0^T \int_K \left( [ \vre \vTe ]_{\rm res} \right)^{\gamma + \beta} \dx \dt \leq c(K) \ep^{2 + \omega}
\ \mbox{for certain}\ \beta, \ \omega > 0 \ \mbox{and any compact}\ K \subset R^3.
\end{equation}

\begin{Remark} \label{RRRH}

Note that estimate (\ref{r12}) is only local in space. It is only at this moment of the proof, where we effectively use
the restriction
\[
0 < \alpha < \frac{4}{3}
\]
imposed by (\ref{HYP}).

\end{Remark}

\section{Acoustic waves}
\label{A}

The present section represents the bulk of the paper. Our goal is to derive dispersive estimates for the acoustic waves. 
As it is probably impossible to derive estimates that would be uniform in the whole frequency spectrum, we combine the perturbative result of 
Metcalfe and Tataru \cite{MetTat1}, with the decomposition technique introduced by Smith and Sogge \cite{SmSo}. 
To this end, we decompose 
the initial data as
\[
\begin{split}
\vr(0, \cdot) &= \vr_0 + \ep ( \vr^{(1)}_{0,\ep} - \vr^{(1)}_0 ) + \ep \vr^{(1)}_0,\\
\vu(0, \cdot) &= \vu_{0,\ep} - \vc{v}_0 + \vc{H}_{\vr_0}[\vc{v}_0] + \vc{v}_0 - \vc{H}_{\vr_0}[\vc{v}_0] = \vu_{0,\ep} - \vc{v}_0 + \vc{V}_0 + \Grad \Phi_0,
\end{split}
\]
where, by virtue of (\ref{i10}), (\ref{i11}), the terms $\vr_0$ and $\vc{V}_0$ represent the zero-th order approximation,
the terms $( \vr^{(1)}_{0,\ep} - \vr^{(1)}_0 )$ and $\vu_{0,\ep} - \vc{v}_0$ vanish in the asymptotic limit, and
\[
\vr^{(1)}_0, \ \Grad \Phi_0, \ \mbox{where}\  \Div \left( \vr_0 \Grad \Phi_0 \right) = \Div ( \vr_0 \vc{v}_0)
\]
create the acoustic waves governed by the system of equations
\begin{eqnarray}
\label{A1}
\ep \partial_t s_\ep + \Div \left[ \vr_0 \Grad \Phi_\ep \right] &=& 0,\\
\label{A2}
\ep \vr_0 \partial_t \Grad \Phi_\ep + \vr_0 \Grad \left[ \frac{p'(\vr_0)}{\vr_0}
s_\ep \right] &=& 0,
\end{eqnarray}
with the initial data
\begin{equation} \label{A3}
s_\ep (0, \cdot) = \vr^{(1)}_0, \ \Grad \Phi_\ep (0, \cdot) = \Grad \Phi_0.
\end{equation}

Note that (\ref{A1}--\ref{A3}) is exactly the same system as obtained in \cite{FeJiNo1}, where the analysis has been performed under a highly simplifying
assumption of $F$ being \emph{compactly supported}. In such a case, the density profile $\vr_0$ given by (\ref{i4}) coincides with the constant $\Ov{\vr}$ outside
a bounded ball and the problem can be treated as a compact perturbation of the standard wave equation generated by ``flat'' Laplacian. In the present setting,
$\vr_0$ shares the asymptotic properties of the long range potential  $F$ specified in (\ref{i4a}--\ref{i4c}), more precisely
\begin{equation} \label{A5}
\vr_0(x) = Q^{-1} \left( F(x) + Q (\Ov{\vr}) \right), \ \mbox{where}\ Q'(r) = \frac{p'(r)}{r} = \gamma r^{\gamma - 2}.
\end{equation}

We introduce the differential operator
\[
\mathcal{A}_{\vr_0}: v \mapsto - \frac{p'(\vr_0)}{\vr_0} \Div \left( \vr_0 \Grad v \right)
\]
that can be interpreted as a non-negative, self-adjoint operator on the weighted space $L^2(R^3)$ space endowed with the scalar product
\[
\left< u; v \right> = \intR{ u v \frac{\vr_0}{p'(\vr_0)}  },
\]
see DeBi\` evre and Pravica \cite{DeBpr1}. Moreover,
as $\vr_0$ is given by (\ref{A5}), where $F$ satisfies (\ref{i4a}--\ref{i4c}), we may infer that
\begin{itemize}
\item
the point spectrum of $\mathcal{A}_{\vr_0}$ is empty, see DeBi\` evre and Pravica \cite[Theorem 2.1(a)]{DeBpr2};
\item
the operator $\mathcal{A}_{\vr_0}$ satisfies the \emph{Limiting absorption principle}, see DeBi\` evre and Pravica \cite[Theorem 1.1]{DeBpr1}.
\end{itemize}

Finally, we may use the same arguments as in  \cite[Section 6.3]{EF80} (see also \cite{DoFeNo2010}), in particular the abstract theorem of Kato \cite{Kato1}, to show the frequency localized decay estimates:
\begin{equation} \label{A4}
\int_{-\infty}^{\infty} \left\| \varphi G(\mathcal{A}_{\vr_0}) \exp\left( \pm {\rm i} \sqrt{\mathcal{A}_{\vr_0}} t \right) [h] \right\|_{L^2(R^3)}^2 \ \dt
\leq c(\varphi, G) \| h \|_{L^2(R^3)}^2
\end{equation}
for any $\varphi \in \DC(R^3)$ and any $G \in \DC(0, \infty)$. Note that $\varphi$ represents a cut-off in the physical space while $G$ a cut-off in
the frequency space. Such an estimate is too weak to be used in combination with the relative energy inequality requiring stronger estimates of
Strichartz's type, cf. \cite{FeJiNo1}. They will be derived in the remaining part of the section.

\subsection{Strichartz estimates by Metcalfe and Tataru}

Metcalfe and Tataru \cite{MetTat1} proved global Strichartz estimates for solutions of the wave equation
\begin{equation} \label{A14}
\partial^2_{t,t} V - \Div \left( \tilde A(x) \Grad V \right) + \tilde{\vc{B}} (x) \cdot \Grad V = Z,\
V(0, \cdot) = V_0, \ \partial_t V(0, \cdot) = V_1
\end{equation}
in the form
\begin{equation} \label{A15}
\| V \|_{L^p(-\infty, \infty; L^q(R^3))} \leq \left( \| \Grad V_0 \|_{L^2(R^3;R^3)} + \| V_1 \|_{L^2(R^3)}
+ \| Z \|_{L^r(-\infty, \infty; L^s(R^3))}   \right),
\end{equation}
whenever
\begin{equation} \label{A16}
\frac{1}{p} + \frac{3}{q} = \frac{1}{2} = \frac{1}{r} + \frac{3}{s} - 2.
\end{equation}
The coefficients $\tilde A(x)$, $\tilde{\vc{B}}(x)$ must be ``asymptotically flat'' in the sense that
\begin{eqnarray} \label{A16a}
\sum_{j \in Z} \sup_{x \in A_j} \left( |x|^2 |\Grad^2 \tilde A(x) | + |x| |\tilde A(x)| + |\tilde A(x) - \Ov{A} | \right) &\leq& \delta,\
\Ov{A} > 0, \\ \label{A16b}
\sum_{j \in Z} \sup_{x \in A_j} \left( |x|^2 |\Grad \tilde{\vc{B}} (x) | + |x| |\tilde{\vc{B}} (x)| \right) &\leq& \delta,
\end{eqnarray}
where $A_j = \{ 2^j \leq |x| \leq 2^{j+1} \}$ are the dyadic regions covering $R^3$, and $\delta > 0$ is sufficiently small, see
Metcalfe and Tataru \cite[Theorem 2]{MetTat1}.

Now, we rewrite the operator $\mathcal{A}_{\vr_0}$ in the form used by Metcalfe and Tataru, specifically,
\[
\begin{split}
\mathcal{A}_{\vr_0}[v] &= - \Div \left( p'(\vr_0) \Grad v \right) +  \vr_0 \Grad \left( \frac{p'(\vr_0)}{\vr_0} \right) \cdot \Grad v
\\
& = - \Div \left( p'(\vr_0) \Grad v \right) + \vr_0 Q''(\vr_0) \Grad \vr_0 \cdot \Grad v,
\end{split}
\]
and set
\begin{equation} \label{A6-}
A(x) = p'(\vr_0)(x), \ \vc{B}(x) = \vr_0 Q''(\vr_0) \Grad \vr_0.
\end{equation}
We easily compute that
\begin{equation} \label{A6}
| \Grad A(x) | + |\vc{B}(x) | \leq c |\Grad F (x) | \leq \frac{c}{(1 + |x|^2)},
\end{equation}
and
\begin{equation} \label{A7}
| \Grad^2 A(x) | + |\Grad \vc{B}(x) | \leq c\left( |\Grad^2 F(x)| +  |\Grad F (x) |^2 \right) \leq \frac{c}{(1 + |x|^3)}.
\end{equation}

Of course, the result of Metcalfe and Tataru \cite{MetTat1} may not be applicable to the operator $\mathcal{A}_{\vr_0}$ as $A$ and $\vc{B}$ may fail
to comply with (\ref{A16a}), (\ref{A16b}) for small $\delta$. Instead, we
consider a modified density profile $\tilde \vr$,
\begin{equation} \label{A8}
\tilde \vr(x) = Q^{-1} \left( H_R(F(x)) + Q (\Ov{\vr}) \right),
\end{equation}
where
\begin{equation} \label{A9}
\begin{split}
H_R \in C^\infty(R), \ H_R(z) &= z \ \mbox{on an open interval containing}\ (-\delta, \delta)\\
|H'_R(z)| &\leq 1 \ \mbox{for all}\ z, \ |H''_R(z)| \leq \frac{c}{R}, \ H_R(z) = 0 \ \mbox{for all}\ z \geq R,\\
\ \mbox{with a constant}\ &c \ \mbox{independent of}\ R > 0.
\end{split}
\end{equation}
Next, we introduce a wave operator $\mathcal{A}_{\tilde \vr}$,
\[
\begin{split}
\mathcal{A}_{\tilde \vr}[v]
= \Div \left( p'(\tilde \vr) \Grad v \right) + \tilde \vr Q''(\tilde \vr) \Grad \tilde \vr \cdot \Grad v,
\end{split}
\]
with
\[
\tilde A(x) = p'(\tilde \vr)(x), \ \tilde{\vc{B}}(x) = \tilde \vr Q''(\tilde \vr) \Grad \tilde \vr.
\]
Observe that, since $H_R$ is linear on a neighborhood of zero and $F(x) \to 0$ as $|x| \to \infty$,
\begin{equation} \label{A10-}
\mathcal{A}_{\tilde \vr} = \mathcal{A}_{\vr_0} \ \mbox{on the set}\ |x| > M , \ M = M(R) \ \mbox{large enough}.
\end{equation}

Our goal is to show that for any $\delta > 0$, there is $R>0$ sufficiently small such that $\tilde{A}$, $\tilde{\vc{B}}$ satisfy
(\ref{A16a}), (\ref{A16b}), in other words, we have the Strichartz estimates (\ref{A15}) for equation (\ref{A14}).
Indeed, as $F > 0$, we have that for any $L > 0$, there is $R = R(L)$ such that
\begin{equation} \label{A10}
\tilde{A}(x) = p'(\Ov{\vr}) \ \mbox{for all} \ |x| < L,\ |\tilde{A}(x) - p'(\Ov{\vr})| \leq H_R(F(x)) \leq \frac{1}{|x|}
\ \mbox{for}\ |x| \geq L.
\end{equation}
By the same token, using (\ref{A6}), we get
\begin{equation} \label{A11}
\begin{split}
\Grad \tilde A, \ \tilde {\vc{B}} & = 0 \ \mbox{for all}\ |x| < L,\\
| \Grad \tilde A(x) | + |\tilde{\vc{B}}(x) | &\leq c |\Grad H_R(F) (x) | \leq c |\Grad F | \leq \frac{c}{|x|^2}
\ \mbox{for}\ |x| \geq L.
\end{split}
\end{equation}
Finally, by virtue of (\ref{A7}), we also have
\begin{equation} \label{A12}
\begin{split}
\Grad^2 \tilde A, \ \Grad \tilde {\vc{B}} & = 0 \ \mbox{for all}\ |x| < L,\\
| \Grad^2 \tilde A(x) | + |\Grad \tilde{\vc{B}}(x) | &\leq c \left( |\Grad^2 H_R(F) (x) | + |\Grad H_R(F)(x) |^2 \right) \\
&\leq c \left[ (1 + H''_R(F(x)) ) |\Grad F(x) |^2  + |\Grad^2 F(x) |^2 \right] \\
&\leq c \left[ \frac{1}{|x|^3} +  H''_R(F(x)) \frac{1}{|x|^4} \right]
\ \mbox{for}\ |x| \geq L.
\end{split}
\end{equation}
On the other hand, however,
\[
H''_R(F) \ne 0 \ \Rightarrow \ 0 \leq F(x) \leq R,\
\]
and, as $F(x) \approx \frac{1}{|x|}$ for $|x| \to \infty$ and by virtue of (\ref{A9}),
\[
|H''_R(F)(x)| \leq \frac{c} {R} \leq c|x|;
\]
whence (\ref{A12}) yields
\begin{equation} \label{A13}
\begin{split}
\Grad \tilde A,\
\Grad^2 \tilde A, \ \Grad \tilde {\vc{B}} & = 0 \ \mbox{for all}\ |x| < L,\\
| \Grad^2 \tilde A(x) | + |\Grad \tilde{\vc{B}}(x) | &\leq
\frac{c}{|x|^3}
\ \mbox{for}\ |x| \geq L.
\end{split}
\end{equation}
Thus if $R = R(\delta) > 0$ is taken small enough, the coefficients $\tilde{A}$, $\tilde{\vc{B}}$ satisfy
(\ref{A16a}), (\ref{A16b}) and the result of
Metcalfe and Tataru \cite{MetTat1} applies yielding (\ref{A14}--\ref{A16}).

\subsection{Decomposition method of Smith and Sogge}

Our final goal is to combine the localized decay estimates (\ref{A4}) with the global ones established in the preceding section to derive
\emph{frequency localized} Strichartz estimates for the operator $\mathcal{A}_{\vr_0}$ in the form
\begin{equation} \label{A18}
\int_{-\infty}^{\infty} \left\| G(\mathcal{A}_{\vr_0}) \exp\left( \pm {\rm i} \sqrt{\mathcal{A}_{\vr_0}} t \right) [h] \right\|_{L^q(R^3)}^p \ \dt
\leq c(G) \| h \|_{L^2(R^3)}^p
\end{equation}
for any any $G \in \DC(0, \infty)$ provided
\[
\frac{1}{p} + \frac{3}{q} = \frac{1}{2}.
\]

\begin{Remark} \label{RA8}

The main difference between (\ref{A4}) and (\ref{A18}) is the absence of the spatial cut-off in the latter. As the functions
$G(-\mathcal{A}_{\vr_0}) \exp\left( \pm {\rm i} \sqrt{\mathcal{A}_{\vr_0}} t \right) [h]$ are smooth, estimate (\ref{A18}) yields uniform decay in
$x$ of spatial derivatives of arbitrary order.

\end{Remark}

To deduce (\ref{A18}), we use the decomposition method by Smith and Sogge \cite{SmSo}.
We consider
\[
U = G(\mathcal{A}_{\vr_0} ) \exp \left( \pm {\rm i} \sqrt{ \mathcal{A}_{\vr_0}  } t \right)[h]
\]
- the solution of the problem
\[
\begin{split}
\partial^2_{t,t} U &- \Div \left( A(x) \Grad U \right) + {\vc{B}} (x) \cdot \Grad U = 0,\\
& U(0, \cdot) = G(\mathcal{A}_{\vr_0} )[h], \ \partial_t U(0, \cdot) = \pm {\rm i} \sqrt{ \mathcal{A}_{\vr_0}  } G(- \mathcal{A}_{\vr_0} )[h],
\end{split}
\]
where $A$ and $\vc{B}$ have been introduced in (\ref{A6-}).

Consider a cut-off function
\[
\chi \in \DC(R), \ \chi(z) = \left\{ \begin{array}{l} 1 \ \mbox{for}\ |z| < M + 1 \\ 0 \ \mbox{for}\ |z| > M+2, \end{array} \right.
\]
where $M$ is from (\ref{A10-}).
Now, we decompose
\[
U = \chi(|x|) U + (1 - \chi (|x|)) U,
\]
where, by virtue of (\ref{A10-}), the function $V = (1 - \chi(|x|)) U$ satisfies
\begin{equation} \label{A17}
\begin{split}
\partial^2_{t,t} V - \Div \left( \tilde A(x) \Grad V \right) &+ \tilde{\vc{B}} (x) \cdot \Grad V = Z,\\
Z &=  \tilde{A} \Grad \chi \cdot \Grad U + \Div ( \tilde A \Grad \chi U ) - \tilde{\vc{B}} \cdot \Grad \chi U \\
V(0, \cdot) &= (1 - \chi) U_0 , \ \partial_t V(0, \cdot) = (1 - \chi) U_1.
\end{split}
\end{equation}

Now, as $U$ satisfies the local estimates (\ref{A4}), the desired decay properties for $\chi U$ follow. As for $V$, we can apply
the Strichartz estimates (\ref{A15}) for $r = 2$, $s = \frac{3}{2}$, estimating $Z$ by means of (\ref{A4}). Note that all components of $Z$ have
compact support in $R^3$ and $U$ as well as $\tilde A$, $\tilde{\vc{B}}$, and $\chi$ are smooth. Thus we have proved
(\ref{A18}).

\subsection{Dispersive estimates for the rescaled equation}

Going back to the original rescaled system (\ref{A1}--\ref{A3}) we first regularize the initial data in the way similar to \cite{EF100}, \cite{FeJiNo1}
taking
\begin{equation} \label{A20}
s_\ep (0, \cdot) = s_{0, \delta} = \frac{\vr_0}{p'(\vr_0)} \left[ \frac{p'(\vr_0)}{\vr_0} \vr^{(1)}_0 \right]_\delta,\
\Phi_\ep (0 \cdot) = \Phi_{0, \delta} = \left[ \Phi_{0} \right]_\delta,
\end{equation}
where the regularization operator $[h]_\delta$ is defined as
\[
[h]_\delta = G_\delta \left(\sqrt{ \mathcal{A}_{\vr_0} } \right)[ \psi_\delta h ],
\]
\[
\psi_\delta \in \DC(R^3), \ 0 \leq \psi_\delta \leq 1, \ \psi_\delta(x) = 1 \ \mbox{for} \ |x| < \frac{1}{\delta},
\ \psi_\delta(x) = 0 \ \mbox{for}\ |x| > \frac{2}{\delta},
\]
\[
\begin{split}
G_\delta &\in \DC(R),\
0 \leq G_\delta \leq 1, \ G_\delta (-z) = G_\delta(z),
\\
G_\delta(z) &= 1 \ \mbox{for}\ z \in \left( - \frac{1}{\delta}, - \delta \right) \cup \left( {\delta},
\frac{1}{\delta}
\right), \\ G_\delta(z) &= 0 \ \mbox{for}\ z \in \left(-\infty, - \frac{2}{\delta} \right) \cup
\left(- \frac{\delta}{2}, \frac{\delta}{2} \right) \cup \left( \frac{2}{ \delta}, \infty \right).
\end{split}
\]

Finally, exactly as in \cite{FeJiNo1}, we may combine the frequency localized estimate
(\ref{A18}), together with the standard energy estimates for the rescaled system (\ref{A1}--\ref{A3}), to
conclude that
\begin{equation} \label{A21}
\begin{split}
\sup_{t \in [0,T]} &\left( \left\| \Phi_\ep (t, \cdot) \right\|_{W^{m,2}(R^3)} + \left\| s_\ep (t, \cdot) \right\|_{W^{m,2}(R^3)} \right) \\
&\leq
c_E(m, \delta) \left( \left\| \Phi_{0,\delta} \right\|_{L^{2}(R^3)} + \left\| s_{0, \delta} \right\|_{L^{2}(R^3)} \right)
\\
\int_0^T &\left( \left\| \Phi_\ep (t, \cdot) \right\|_{W^{m,\infty}(R^3)} + \left\| s_\ep (t, \cdot) \right\|_{W^{m,\infty}(R^3)} \right) \dt\\
&\leq c_D(\ep, m, \delta) \left( \left\| \Grad \Phi_{0,\delta} \right\|_{L^{2}(R^3;R^3)} + \left\| s_{0, \delta} \right\|_{L^{2}(R^3)} \right),
\end{split}
\end{equation}
for any $m$ and $\delta > 0$, where
\[
c_D(\ep, m, \delta) \to 0  \ \mbox{as}\ \ep \to 0 \ \mbox{for}\ m, \ \delta \ \mbox{fixed}.
\]

\section{Conclusion, proof of the main result}
\label{C}

The proof of Theorem \ref{Tm1} will be concluded by taking the limit solution in place of the test functions in the relative energy inequality (\ref{r1}).
To this end, it is convenient to rewrite (\ref{i1}--\ref{i3}) in terms of the ``temperature''
\[
\mathcal{T} = \frac{\vr_0}{\mathcal{R}},
\]
namely
\begin{eqnarray}
\Div (\vr_0 \vc{V}) &=& 0 \label{C1} \\
\partial_t \vc{V} + \vc{V} \cdot \Grad \vc{V} + \Grad \Pi &=& - \mathcal{T} \Grad F \label{C2} \\
\partial_t \mathcal{T} + \vc{V} \cdot \Grad \mathcal{T}&=& 0 \label{C3},
\end{eqnarray}
with the initial data
\begin{equation} \label{C4}
\vc{V}(0, \cdot) = \vc{H}_{\vr_0}[\vc{v}_0], \ \mathcal{T}(0, \cdot) = \Theta^{(2)}_0.
\end{equation}

\subsection{Relative energy ansatz}

We consider
\[
r = r_\ep = \vr_0 + \ep s_\ep, \ \vc{U} = \vc{U}_\ep = \vc{V} + \Grad \Phi_\ep,
\]
where $[s_{\ep}, \Phi_\ep]$ solve the acoustic wave equation (\ref{A1}), (\ref{A2}), and $\vc{V}$ is the velocity in the limit problem (\ref{C1}--\ref{C4}),
as test functions in the relative energy inequality (\ref{r1}).
Note that the initial data can be replaced by their regularizations (\ref{A20}) as
\[
s_{0,\delta} \to \vr^{(1)}_0 \ \mbox{in}\ L^2(R^3) , \ \Grad \Phi_{0, \delta} \to \Grad \Phi_0
\ \mbox{in}\ L^2(R^3; R^3) \ \mbox{as}\ \delta \to 0.
\]
With $\delta > 0$ fixed,
our ultimate goal will be to show that all terms on the right-hand side of the relative energy inequality (\ref{r1}) can be either ``absorbed'' by the left-hand
side by means of a Gronwall type argument or vanish in the asymptotic limit for $\ep \to 0$. This will be done in several steps.

\subsubsection{Step 1 - rewriting the pressure terms}

To begin, observe that
\[
\partial_t r_\ep + \Div ( r_\ep \vc{U}_\ep ) = \ep \partial_t s_\ep + \Div \left( r_\ep \vc{V} + r_\ep \Grad \Phi_\ep \right) =
\ep \Div (s_\ep \vc{U}_\ep).
\]
Using this fact we get, after a straightforward manipulation,
\[
\begin{split}
&(r_\ep - \vre \vTe) \partial_t H'(r_\ep) + \Grad H'(r_\ep ) \cdot (r_\ep \vc{U}_\ep - \vre \vTe \vue ) - \Div \vc{U}_\ep \Big( (\vre \vTe)^\gamma - r^\gamma_\ep \Big)
\\
&= - \Div \vc{U}_\ep \Big[ p(\vre \vTe) - p'(r_\ep) (\vre \vTe - r_\ep ) - p(r_\ep) \Big] + \ep (r_\ep - \vre \vTe)H''(r_\ep) \Div (s_\ep \vc{U}_\ep) \\
&+ \vre \vTe \Grad H'(r_\ep) \cdot (\vc{U}_\ep - \vue ).
\end{split}
\]
Furthermore,
\[
\begin{split}
&\vre \vTe \Grad H'(r_\ep) \cdot (\vc{U}_\ep - \vue ) + \vre \Grad F \cdot (\vue - \vc{U}_\ep)\\
&= \vre \vTe \Grad \left[ H'(r_\ep) - H''(\vr_e) (r_\ep - \vr_0) - H'(\vr_0) \right] \cdot (\vc{U}_\ep - \vue) \\
&+
\vre \vTe \Grad \left[ H''(\vr_0) (r_\ep - \vr_0) + H'(\vr_0) \right] \cdot (\vc{U}_\ep - \vue)  + \vre \Grad H'(\vr_0) \cdot (\vue - \vc{U}_\ep)\\
&= \vre \vTe \Grad \left[ H'(r_\ep) - H''(\vr_0) (r_\ep - \vr_0) - H'(\vr_0) \right] \cdot (\vc{U}_\ep - \vue) \\
&+ \vre \vTe \Grad \left[ H''(\vr_0) (r_\ep - \vr_0) \right] \cdot (\vc{U}_\ep - \vue) + \vre \left(1 - \vTe
\right) \Grad H'(\vr_0) \cdot (\vue - \vc{U}_\ep).
\end{split}
\]
Finally, using (\ref{A2}), we conclude
\[
\vre \vTe \Grad \left[ H''(\vr_0) (r_\ep - \vr_0) \right] \cdot (\vc{U}_\ep - \vue) = -\ep^2 \vre \vTe \partial_t \Grad \Phi_\ep \cdot (\vc{U}_\ep - \vue).
\]

Summarizing the previous relations we may rewrite (\ref{r1}) in the form
\begin{equation} \label{C5}
\begin{split}
&\left[ \mathcal{E} \left( \vre, \vTe, \vue \Big| r_\ep , \vc{U}_\ep \right) \right]_{t = 0}^{t = \tau}
+ \ep^\alpha \int_{0}^{\tau}
\intR{ \mathbb{S} (\Grad (\vue-\vc{U}_\ep)) : \Grad (\vue - \vc{U}_\ep)} \dt
\\
&\leq
\int_{0}^{\tau}
\intR{  \vre \left( \partial_t \vc{V} + \vue \cdot \Grad \vc{U}_\ep \right) \cdot \left( \vc{U}_\ep - \vue \right)
+ \ep^\alpha \mathbb {S} (\Grad \vc{U}_\ep) : \Grad (\vc{U}_\ep - \vue ) } \dt
\\
& - \frac{1}{\ep^2} \int_0^\tau \intR{ \Div \vc{U}_\ep \Big[ p(\vre \vTe) - p'(r_\ep) (\vre \vTe - r_\ep ) - p(r_\ep) \Big] }
\dt
\\
&+ \frac{1}{\ep^2} \int_0^\tau \intR{ \vre \vTe \Grad \left[ H'(r_\ep) - H''(\vr_0) (r_\ep - \vr_0) - H'(\vr_0) \right] \cdot (\vc{U}_\ep - \vue) } \ \dt
\\
& + \frac{1}{\ep} \int_0^\tau \intR{ (r_\ep - \vre \vTe)H''(r_\ep) \Div (s_\ep \vc{U}_\ep) } \dt\\
&+\int_0^\tau \intR{ \left[ \vre (1 - \vTe) \partial_t \Grad \Phi_\ep \cdot (\vc{U}_\ep - \vue) + \frac{\vre}{\ep^2} \left(1 - \vTe
\right) \Grad H'(\vr_0) \cdot (\vue - \vc{U}_\ep) \right] } \dt
\end{split}
\end{equation}

\subsubsection{Step 2 - the diffusion term}

It follows from estimate (\ref{r8}) that
\[
\ep^\alpha \int_{0}^{\tau}
\intR{\mathbb {S} (\Grad \vc{U}_\ep) : \Grad (\vc{U}_\ep - \vue ) } \dt \to 0  \ \mbox{as}\ \ep \to 0.
\]

\subsubsection{Step 3 - the velocity term}

We write
\[
\begin{split}
\int_{0}^{\tau}&
\intR{  \vre \left( \partial_t \vc{V} + \vue \cdot \Grad \vc{U}_\ep \right) \cdot \left( \vc{U}_\ep - \vue \right)} \dt \\
& = \int_{0}^{\tau} \intR{ \vre \left(\vue - \vc{U}_\ep \right) \cdot \Grad \vc{U}_\ep \cdot \left( \vc{U}_\ep - \vue \right)} \dt\\
& + \int_{0}^{\tau} \intR{ \vre \left( \vc{V} \cdot \Grad^2 \Phi_\ep + \Grad \Phi_\ep \cdot \Grad \vc{U}_\ep       \right) \cdot \left( \vc{U}_\ep - \vue \right)} \dt
\\
& +\int_{0}^{\tau}
\intR{  \vre \left( \partial_t \vc{V} + \vc{V} \cdot \Grad \vc{V} \right) \cdot \left( \vc{U}_\ep - \vue \right)} \dt,
\end{split}
\]
where the first integral on the right-hand side can be controlled by a similar term in $\mathcal{E}$.

Next, the second integral can be handled as
\[
\begin{split}
&\left| \intR{ \vre \left( \vc{V} \cdot \Grad^2 \Phi_\ep + \Grad \Phi_\ep \cdot \Grad \vc{U}_\ep       \right) \cdot \left( \vc{U}_\ep - \vue \right)}  \right| \\ &\leq \| \Grad \Phi_\ep \|_{W^{1,\infty}(R^3; R^3)}^2 \left( \| \vc{V} \|_{L^\infty(R^3;R^3)} + \| \Grad \vc{U}_\ep \|_{L^\infty(R^3; R^{3 \times 3})} \right)^2
+ c \mathcal{E},
\end{split}
\]
where the first term vanishes in the asymptotic limit thanks to the dispersive estimates (\ref{A21}).

The last term reads
\[
\begin{split}
&\int_{0}^{\tau}
\intR{  \vre \left( \partial_t \vc{V} + \vc{V} \cdot \Grad \vc{V} \right) \cdot \left( \vc{U}_\ep - \vue \right)} \dt\\
&- \int_{0}^{\tau} \intR{  \vre \Grad \Pi \cdot \left( \vc{U}_\ep - \vue \right)} \dt - \int_{0}^{\tau} \intR{  \vre \mathcal{T} \Grad F \cdot \left( \vc{U}_\ep - \vue \right)} \dt,
\end{split}
\]
where
\[
\begin{split}
&\int_{0}^{\tau} \intR{  \vre \Grad \Pi \cdot \left( \vc{U}_\ep - \vue \right)} \dt \\
&= - \int_0^\tau \intO{ \vre \vue \cdot \Grad \Pi } \ \dt + \int_0^\tau \intO{ \vre \Big( \vc{V} + \Grad \Phi_\ep \Big)
\cdot \Grad \Pi } \ \dt.
\end{split}
\]
In view of (\ref{r2}), (\ref{r11}),
\[
\vre \vue \to \Ov{ \vr \vu } \ \mbox{weakly-(*) in}\ L^\infty(0,T; (L^2 + L^{2 \gamma/ (\gamma + 1)}(R^3;R^3)),
\]
where, as can be deduced from the equation of continuity and (\ref{r11}),
\[
\Div \Ov{ \vr \vu } = 0.
\]
Consequently,
\[
\int_0^\tau \intO{ \vre \vue \cdot \Grad \Pi } \ \dt \to 0 \ \mbox{as}\ \ep \to 0.
\]
Similarly, we may combine the dispersive estimates (\ref{A21}) with the anelastic constraint (\ref{i1}) to deduce that
\[
\int_0^\tau \intO{ \vre \Big( \vc{V} + \Grad \Phi_\ep \Big)
\cdot \Grad \Pi } \ \dt \to 0 \ \mbox{as}\ \ep \to 0.
\]

Thus (\ref{C5}) can be reduced to
\begin{equation} \label{C6}
\begin{split}
&\left[ \mathcal{E} \left( \vre, \vTe, \vue \Big| r_\ep , \vc{U}_\ep \right) \right]_{t = 0}^{t = \tau}
\\
&\leq \frac{1}{\ep^2} \int_0^\tau \intR{ \vre \vTe \Grad \left[ H'(r_\ep) - H''(\vr_0) (r_\ep - \vr_0) - H'(\vr_0) \right] \cdot (\vc{U}_\ep - \vue) } \ \dt
\\
& - \frac{1}{\ep^2} \int_0^\tau \intR{ \Div \vc{U}_\ep \Big[ p(\vre \vTe) - p'(r_\ep) (\vre \vTe - r_\ep ) - p(r_\ep) \Big] }
\dt
\\
& + \frac{1}{\ep} \int_0^\tau \intR{ (r_\ep - \vre \vTe)H''(r_\ep) \Div (s_\ep \vc{U}_\ep) } \dt\\
&+\int_0^\tau \intR{ \vre (1 - \vTe) \partial_t \Grad \Phi_\ep \cdot (\vc{U}_\ep - \vue) } \dt\\
& + \int_0^\tau \intR{ \vre  \left( \frac{1 - \vTe}{\ep^2} - \mathcal{T} \
\right) \Grad F \cdot (\vue - \vc{U}_\ep)  } \dt \\
&+ c \int_0^\tau \mathcal{E} \left( \vre, \vTe, \vue \Big| r_\ep , \vc{U}_\ep \right) \ \dt + o(\ep),
\end{split}
\end{equation}
where
\[
o(\ep) \to 0 \ \mbox{as}\ \ep \to 0.
\]

\subsubsection{Step 4 - controlling pressure, I}

We compute
\[
\begin{split}
\Grad \Big( H'(r_\ep)&-H''(\vr_0) (r_\ep -
\vr_0) - H'(\vr_0) \Big) \\
&=\ep \Big(H''(r_\ep)-H''(\vr_0)\Big)\Grad
s_\ep+
\Big(P''(r_\ep)-P''(\vr_0)-P'''(\vr_0)(r_\ep- \vr_0)\Big)\Grad \vr_0,
\end{split}
\]
therefore we may use Taylor's formula to deduce
\[
\left| \Grad \Big( H'(r_\ep) -H''(\vr_0) (r_\ep -
\vr_0) - H'(\vr_0) \Big) \right| \le \ep^2 c \Big(s_\ep|\Grad s_\ep|+s_\ep^2\Big).
\]
Thus we conclude
\[
\begin{split}
\frac{1}{\ep^2}  & \left| \int_0^\tau \intR{ \vre \vTe \Grad \left[ H'(r_\ep) - H''(\vr_0) (r_\ep - \vr_0) - H'(\vr_0) \right] \cdot (\vc{U}_\ep - \vue) } \ \dt
\right|\\
& \leq c \int_0^\tau \intR{ \vre \Big(s_\ep|\Grad s_\ep|+s_\ep^2\Big) |\vue - \vc{U}_\ep |    }\dt,
\end{split}
\]
where, by virtue of the dispersive estimates (\ref{A21}), the last integral is dominated by
\[
o(\ep) + c \int_0^\tau \mathcal{E} \left( \vre, \vTe, \vue \Big| r_\ep , \vc{U}_\ep \right) \ \dt
\]
Relation (\ref{C6}) reduces to
\begin{equation} \label{C7}
\begin{split}
&\left[ \mathcal{E} \left( \vre, \vTe, \vue \Big| r_\ep , \vc{U}_\ep \right) \right]_{t = 0}^{t = \tau}
\\
& - \frac{1}{\ep^2} \int_0^\tau \intR{ \Div \vc{U}_\ep \Big[ p(\vre \vTe) - p'(r_\ep) (\vre \vTe - r_\ep ) - p(r_\ep) \Big] }
\dt
\\
& + \frac{1}{\ep} \int_0^\tau \intR{ (r_\ep - \vre \vTe)H''(r_\ep) \Div (s_\ep \vc{U}_\ep) } \dt\\
&+\int_0^\tau \intR{ \vre (1 - \vTe) \partial_t \Grad \Phi_\ep \cdot (\vc{U}_\ep - \vue) } \dt\\
& + \int_0^\tau \intR{ \vre  \left( \frac{1 - \vTe}{\ep^2} - \mathcal{T} \
\right) \Grad F \cdot (\vue - \vc{U}_\ep)  } \dt \\
&+ c \int_0^\tau \mathcal{E} \left( \vre, \vTe, \vue \Big| r_\ep , \vc{U}_\ep \right) \ \dt + o(\ep),
\end{split}
\end{equation}

\subsubsection{Step 5 - controlling pressure, II}

Probably the most delicate step is to handle the term
\begin{equation} \label{C8}
\begin{split}
\frac{1}{\ep^2} &\int_0^\tau \intR{ \Div \vc{U}_\ep \Big[ p(\vre \vTe) - p'(r_\ep) (\vre \vTe - r_\ep ) - p(r_\ep) \Big] }
\dt\\
&=\frac{1}{\ep^2} \int_0^\tau \intR{ \Big( \Div \vc{V} +        \Delta_x \Phi_\ep \Big) \Big[ p(\vre \vTe) - p'(\vr_0) (\vre \vTe - \vr_0 ) - p(\vr_0) \Big] }
\dt
\\
&+\frac{1}{\ep^2} \int_0^\tau \intR{ \Big( \Div \vc{V} +        \Delta_x \Phi_\ep \Big) \Big[ p(\vr_0) - p(\vr_0 + \ep s_\ep) + \ep
p'(\vr_0 + \ep s_\ep ) s_\ep  \Big] }
\dt
\\
&+\frac{1}{\ep^2} \int_0^\tau \intR{ \Big( \Div \vc{V} +  \Delta_x \Phi_\ep \Big) \Big[ p'(\vr_0) - p'(\vr_0 + \ep s_\ep)  \Big]
(\vre \vTe - \vr_0) }
\dt,
\end{split}
\end{equation}
where
\[
\left| \frac{ p(\vr_0) - p(\vr_0 + \ep s_\ep) + \ep
p'(\vr_0 + \ep s_\ep ) s_\ep }{\ep^2} \right| \leq c |s_\ep |^2,
\]
and, similarly,
\[
\left|
\frac{ p'(\vr_0) - p'(\vr_0 + \ep s_\ep)}{\ep} \right| \left|
\frac{ (\vre \vTe - \vr_0) }{\ep} \right| \leq c |s_\ep | \left|
\frac{ (\vre \vTe - \vr_0) }{\ep} \right|.
\]
Thus the last two integrals in (\ref{C8}) can be controlled by means of (\ref{r6}) and the dispersive estimates (\ref{A21}).

Next, observing that
$\Delta_x \Phi_\ep$ vanishes for $\ep \to 0$ because  of (\ref{A21}) we have
\[
\frac 1{\ep^2}\int_0^\tau\intR{ \Delta_x \Phi_\ep \Big( p(\vre \vTe) - p'(\vr_0)
(\vre \vTe - \vr_0) - p(\vr_0) \Big) } \
\dt\to 0\;\mbox{in $L^1(0,T)$ as $\ep\to 0$}.
\]

The next step is to handle
\[
\begin{split}
\frac {1}{\ep^2} &\left| \int_0^\tau\intR{ \Div \vc{V} \Big( p(\vre \vTe) - p'(\vr_0)
(\vre \vTe - \vr_0) - p(\vr_0) \Big)  } \ \dt \right|\\
\\
&\le \frac 1{\ep^2} \int_0^\tau\intR{ \Big|\Big[ p(\vre \vTe) - p'(
\vr_0) (\vre \vTe - \vr_0) - p(\vr_0) \Big]_{\rm ess}\Big|
|\Div \vc{V} | } \ \dt \\
&
\frac 1{\ep^2} \int_0^\tau\intR{ \Big|\Big[ p(\vre \vTe) - p'(
\vr_0) (\vre \vTe - \vr_0) - p(\vr_0) \Big]_{\rm res}\Big|
|\Div \vc{V} | } \ \dt.
\end{split}
\]
As the limit velocity field $\vc{V}$ obeys the anelastic constraint (\ref{i1}), we get
\[
\Div \vc{V} = - \frac{\Grad \vr_0}{\vr_0} \cdot \vc{V}.
\]
Since $\vr_0$ is given by (\ref{i4}), we may find a compact set $K \subset R^3$ such that $\Div \vc{V}$ is small in $R^3 \setminus K$.
Thus we may use the \emph{local} pressure estimates established in (\ref{r12}) to obtain
\[
\begin{split}
\frac {1}{\ep^2} \int_0^\tau\int_K \Big|\Big[ p(\vre \vTe) - p'(
\vr_0) (\vre \vTe - \vr_0) - p(\vr_0) \Big]_{\rm res}\Big|
|\Div \vc{V} | \dx \dt \to 0
\ \mbox{as}\ \ep \to 0,
\end{split}
\]
while the integral over the complement $R^3 \setminus K$ is small because of (\ref{r7}).

\begin{Remark} \label{CR1}

This is the only moment in the proof where the pressure estimates are needed.

\end{Remark}

Finally,
\[
\begin{split}
&\left| \int_0^\tau \intR{\left[ \frac{ p(\vre \vTe) - p'(\vr_0)
(\vre \vTe - \vr_0) - p(\vr_0)}{\ep^2} \right]_{\rm ess} |\Div \vc{V} | } \dt \right|
\\
&
\leq c\int_0^\tau\intR{ \left| \left[
\frac{\vre \vTe - \vr_0}{\ep} \right]_{\rm ess} \right|^2 |\Div \vc{V}| } \dt
\\
&\leq
c\int_0^\tau\intR{ \left( \left[
\frac{\vre \vTe - \ep s_\ep - \vr_0}{\ep} \right]_{\rm ess}^2 + |s_\ep|^2  \right) |\Div \vc{V}| } \dt
\\
&\leq c \int_0^\tau\mathcal{E}\left( \vre, \vTe, \vue \Big| r_\ep.
\vc{U}_\ep \right) + \int_0^\tau\intR{|s_\ep|^2|\Div \vc{V}|}{\rm
d}t.
\end{split}
\]

\subsubsection{Step 6 - conclusion}

To conclude, we write
\[
\begin{split}
&\frac{1}{\ep} \left| \int_0^\tau \intR{ (r_\ep - \vre \vTe)H''(r_\ep) \Div (s_\ep \vc{U}_\ep) } \dt \right| \\
&\leq c\int_0^\tau \int_{R^3} \left(\left| \frac{\vre \vTe - \vr_0 }{\ep} \right| + |s_\ep| \right)(|s_\ep| + |\Grad s_\ep |)
( |\vc{V}| + |\Div \vc{V}| + |\Grad \Phi_\ep| + |\Delta_x \Phi_\ep |  ) \dx \dt
\end{split}
\]
where again the right-hand side vanishes for $\ep \to 0$ due to the dispersive estimates (\ref{A21}).

Finally, the integral
\[
\int_0^\tau \intR{ \vre (1 - \vTe) \partial_t \Grad \Phi_\ep \cdot (\vc{U}_\ep - \vue) } \dt
\]
is small because of (\ref{r5}), and
\[
\int_0^\tau \intR{ \vre  \left( \frac{1 - \vTe}{\ep^2} - \mathcal{T} \
\right) \Grad F \cdot (\vue - \vc{U}_\ep)  } \dt
\]
is controlled by the initial data satisfying
(\ref{i12}).

Thus we may infer that
\[
\begin{split}
\left[ \mathcal{E} \left( \vre, \vTe, \vue \Big| r_\ep , \vc{U}_\ep \right) \right]_{t = 0}^{t = \tau}
\leq  c \int_0^\tau \mathcal{E} \left( \vre, \vTe, \vue \Big| r_\ep , \vc{U}_\ep \right) \ \dt + o(\ep),\\
o(\ep) \to 0 \ \mbox{as}\ \ep \to 0
\end{split}
\]
and use Gronwall's lemma to complete the proof of Theorem \ref{Tm1}.

\def\cprime{$'$} \def\ocirc#1{\ifmmode\setbox0=\hbox{$#1$}\dimen0=\ht0
  \advance\dimen0 by1pt\rlap{\hbox to\wd0{\hss\raise\dimen0
  \hbox{\hskip.2em$\scriptscriptstyle\circ$}\hss}}#1\else {\accent"17 #1}\fi}

\S\S\S\S\S\S\S\S\S\S\S\S\S

\end{document}